\documentclass{amsart}
\usepackage{amsfonts}
\usepackage{enumerate}
\usepackage{latexsym}
\usepackage{amsmath}

\theoremstyle{definition}
\newtheorem{defi}{Definition}[section]
\newtheorem{remark}{Remark}

\theoremstyle{plane}
\newtheorem{theo}{Theorem} 
\newtheorem{lem}{Lemma}[section]
{}

\newtheorem{cor}{Corollary}[section]

\newcommand{\R}{{\mathbb R}}
\newcommand{\Z}{{\mathbb Z}}

\newcommand{\ha}{{\rm Harm}}
\newcommand{\ho}{{\rm Hom}}
\newcommand{\Pol}{{\rm Pol}}
\newcommand{\Ker}{{\rm Ker}}

\title{On Euclidean designs and the potential energy}
\author{Tsuyoshi Miezaki}
\address{Division of Mathematics, 
Graduate School of Information Sciences, 
Tohoku University, 
6-3-09 Aramaki-Aza-Aoba, Aoba-ku, Sendai 980-8579, Japan}
\email{miezaki@math.is.tohoku.ac.jp}
\author{Makoto Tagami}
\address{Mathematical Institute, Tohoku University,
6-3, Aoba, Aramaki, Aoba-ku, Sendai 980-8578, Japan}
\email{tagami@math.tohoku.ac.jp}
\begin{document}
\maketitle
\begin{abstract}
We study Euclidean designs from the viewpoint of the potential energy. For a finite set in Euclidean space, We formulate a linear programming bound for the potential energy by applying harmonic analysis on a sphere. We also introduce the concept of strong Euclidean designs from the viewpoint of the linear programming bound, and we give a Fisher type inequality for strong Euclidean designs. A finite set on Euclidean space is called a Euclidean $a$-code if any distinct two points in the set are separated at least by $a$. As a corollary of the linear programming bound, we give a method to determine an upper bound on the cardinalities of Euclidean $a$-codes on concentric spheres of given radii. Similarly we also give a method to determine a lower bound on the cardinalities of Euclidean $t$-designs as an analogue of the linear programming bound. 
\end{abstract}
\section{Introduction}
The concept of Euclidean designs is well-known as a natural generalization of spherical designs to Euclidean space. The purpose of this paper is firstly to characterize Euclidean designs in terms of the potential energy, secondly to extend the linear programming bounds on a sphere to configurations of points in Euclidean space, and thirdly to introduce the concept of strong Euclidean designs which seems to be natural from the view point of the linear programming bound.

Let $t$ be a natural number, $S^{d-1}$ the $(d-1)$-dimensional unit sphere centered at the origin, and $P_l(\R^d)$ the vector space of polynomials of degree at most $l$ in $d$ variables over $\R$. A finite nonempty subset $X$ on $S^{d-1}$ is called a \textit{spherical $t$-design} if, for any $f(x)\in P_t(\R^d)$, the following equality holds:
\[\frac{1}{|S^{d-1}|}\int_{S^{d-1}}f(x)d\sigma(x)=\frac{1}{|X|} \sum_{x\in X} f(x).\]
Here $\sigma$ is an $O(\R^d)$-invariant measure on $S^{d-1}$ and $|S^{d-1}|$ denotes the surface volume of the sphere $S^{d-1}$. The concept of spherical designs was defined by Delsarte-Goethals-Seidel \cite{DGS}. A spherical $t$-design means to be a good configuration of points on the sphere so that the average value of the integral of any polynomial of degree up to $t$ on the sphere is replaced by the average value at a finite set on the sphere. A finite non-empty subset $X$ on $S^{d-1}(r)$, the sphere of radius $r$  centered at the origin, is also called a spherical $t$-design if $\frac{1}{r}X$ is a spherical $t$-design on the unit sphere $S^{d-1}$.

Let $\triangle$ be the Laplacian, that is, 
$\triangle=\frac{\partial^2}{\partial x_1^2}+\cdots+\frac{\partial^2}{\partial x_d^2}$, and set
\begin{align*}
\ho_i(\R^d)&=\{f(x)\in \R[x_1,\ldots,x_d]\mid f(x)\;\mbox{is homogeneous of 
degree $i$}\}, \\
\ha_i(\R^d)&=\{f(x) \in \ho_i(\R^d)\mid \triangle (f(x))= 0\}. 
\end{align*}
The following is an equivalent condition for $X$ to be a spherical $t$-design.
\begin{lem}[cf. \cite{DGS}]\label{sdec}
A finite non-empty subset $X$ on $S^{d-1}$ is a spherical $t$-design 
if and only if the following condition holds:
\[\sum_{x\in X}\phi(x)=0 \;(\phi \in \ha_j(\R^d),\;1\leq j \leq t).\]
\end{lem}
Let $\Z_{\ge 0}$, $\Z_{>0}$ and $\R_{>0}$ denote the set of non-negative integers, positive integers and positive real numbers, respectively. A spherical $t$-design is closely related to the following Sidelnikov inequality (refer to \cite{sidelnikov}): for a finite subset $X$ on $S^{d-1}$ and any $l \in \Z_{\ge 0}$, it holds that
\begin{equation}\label{sid}
\frac{1}{|X|^2}\sum_{x,y\in X} (x\cdot y)^l \ge A_l := 
\begin{cases} 
\displaystyle \frac{(l-1)!!(d-2)!!}{(d+l-2)!!} & \mbox{if $l\equiv 0 \pmod{2}$},\\
\displaystyle 0 & \mbox{otherwise}. 
\end{cases}
\end{equation}
Here $x\cdot y$ is the standard inner product, and for $x\in \R$, $[x]$ denotes the largest integer not exceeding $x$, and we define $l!!$ to be $\prod_{i=0}^{[\frac{l}{2}]} (l-2i)$ if $l\in\Z_{> 0}$ and $1$ if $l\le 0$. It is well-known that $X$ is a spherical $t$-design if and only if for $0\le l \le t$, equality holds in (\ref{sid}). 

Consider a function $f:(0,4]\rightarrow \R$. Then the \textit{potential energy} of $X$ for $f$ is defined to be $PE_f(X):=\sum_{x\not=y \in X} f(||x-y||^2)$ (refer to \cite {yudin,cohn}). From Sidelnikov's inequality and the equivalent condition, we see that, for $f_t(r)=(4-r)^t$, a spherical $t$-design $X$ minimizes the potential energy $PE_{f_t}(X)$ among all configurations of $|X|$ points on $S^{d-1}$.

The concept of spherical designs was naturally generalized to Euclidean space 
by Neumaier-Seidel \cite{Neu}. Let $X$  be a finite subset in $\R^d$ and suppose $O \not\in X$ (for simplicity, we always suppose this in this paper). We set $RX:=\{||x|| \mid x \in X\}=\{r_1,r_2,\ldots,r_p\}$, $S_i=r_iS^{d-1}$, $RS:=S_1\cup S_2\cup\cdots\cup S_p$ and $X_{i}=X_{(r_i)}=X \cap S_i$. $\sigma_i$ denotes an $O(\R^d)$-invariant measure on $S_i$. Also let $w:X \rightarrow \R_{>0}$ and we put $w(X_i):=\sum_{x\in X_i}w(x)$.
\begin{defi}[Neumaier-Seidel \cite{Neu}]
Under the above notation, $(X,w)$ is a \textit{Euclidean $t$-design} if 
for any $f(x)\in P_t(\R^d)$, the following equality holds:
\[\sum_{i=1}^p \frac{w(X_i)}{|S_i|}\int_{S_i}f(x)d\sigma_i(x)=\sum_{x\in X} w(x)f(x).\]
\end{defi}
In Euclidean space, the following generalized Sidelnikov inequality holds:
\begin{theo}[Neumaier-Seidel \cite{Neu}]\label{gs}
Let $X \subset \R^d\; (|X|< \infty)$, $w:X\rightarrow \R$ 
and $l \in\R_{\geq 0}$. 
Then the following inequality holds: 
\begin{equation}\label{eucsei}
\sum_{x,y\in X} w(x)w(y)(x\cdot y)^l \ge 
A_l\left(\sum_{x \in X}w(x)||x||^l\right)^2.
\end{equation}
\end{theo}
The actual generalized Sidelnikov inequality given by Neumaier-Seidel \cite{Neu} is a more general one, but in this paper we write the inequality in the form 
of Theorem~\ref{gs} for the purpose of viewing a relation to Euclidean designs.

Let $w':X \longrightarrow \R_{>0}$ and set $w(x)=w'(x)||x||^{2j}$ in Theorem \ref{gs}. Then
\begin{equation}\label{sid2}
\sum_{x,y\in X} w'(x)w'(y)(||x||||y||)^{2j}(x\cdot y)^l \ge 
A_l\left(\sum_{x\in X}w'(x)||x||^{2j+l}\right)^2. 
\end{equation}
$(X,w')$ is a Euclidean $t$-design if and only if equality holds in (\ref{sid2}) for all $j\ge0$, $l\ge0$ with $2j+l\le t$
(cf. Lemma \ref{eujo}, the proof of Theorem \ref{gs} in Section \ref{linearse}).

For $x \in \R^d\setminus \{O\}$, we set $x':={x}/||x||$. 
We define the potential energy of a finite set in Euclidean space 
which is not necessarily restricted on the sphere: 
\begin{defi}\label{pe}
Let $X\subset \R^d\setminus \{O\}\; (|X|< \infty)$, $w:X\rightarrow \R$, and $f:\R_{>0}^2\times [-1,1] \rightarrow \R$.  
 Then the \textit{potential energy} of $(X,w)$ for $f$ is defined by 
\[PE_f(X,w)=\sum_{x,y \in X, x\neq y} w(x)w(y)f(||x||,||y||,x'\cdot y').\] 
\end{defi}
In particular, for the case when $f(r,s,t)=(rs)^l t^l$, 
the generalized Sidelnikov inequality gives a lower bound on $PE_f(X,w)$. 

Let $Q_i(t)=Q_i^{(d)}(t)$ be the Gegenbauer polynomial of degree $i$ 
corresponding to the sphere $S^{d-1}$, namely, 
$\{Q_i(t)\}$ are the orthogonal polynomials on the interval $[-1,1]$ 
with respect to the weight function $(1-t^2)^{\frac{d-3}{2}}$. 
In this paper, they are nomalized as $Q_i(1)=\dim \ha_i(\R^d)$. 

In \cite{yudin}, Yudin gave a linear programming bound 
on the potential energy 
using harmonic analysis on the sphere. 
Applying his method, we obtain a lower bound on the potential energy of 
a finite set on concentric spheres: 
\begin{theo}\label{linear pro}
Let $X\subset \R^d\setminus \{O\}\; (|X|< \infty)$ and $w:X\rightarrow \R$, 
and let $f:\R_{>0}^2\times [-1,1] \rightarrow \R$, 
$g_0:\R_{>0}^2\rightarrow \R$, and $g_{ij}:\R_{>0}\rightarrow \R\;(i,j\ge 1)$. 
If $h(r,s,t)=g_0(r,s)+\sum_{i,j\ge 1} g_{ij}(r) g_{ij}(s)Q_j(t)$ satisfies the condition that 
\[f(r,s,t)\ge h(r,s,t),\;(r,s \in RX,\ t \in [-1,1]),\]
then the following inequality holds: 
\begin{equation}\label{linear}
PE_{f}(X,w)\ge \sum_{x,y \in X} w(x)w(y) g_{0}(||x||,||y||)-\sum_{x\in X}w(x)^2h(||x||,||x||,1).
\end{equation}
Moreover equality holds in $(\ref{linear})$ if and only if 
for any $x,\,y \in X\, (x\neq y)$, 
$f(||x||,||y||,x'\cdot y')=h(||x||,||y||,x'\cdot y')$, and 
\begin{equation}\label{to}
\sum_{x,y \in X} w(x)w(y)g_{ij}(||x||)g_{ij}(||y||)Q_j(x'\cdot y')=0 \; (\forall i\ge 0,\forall j\ge 1).
\end{equation}
\end{theo}
We give a proof of Theorem \ref{linear pro} in Section \ref{linearse}.

Suppose that $RX$ and $|X_{(r)}|$ are given and suppose that  $g_0(r,s)$ is a polynomial and 
each $g_{ij}(r)$ is a monomial. Then seeking $g_0$ and $g_{ij}$ 
which maximize the lower bound (\ref{linear}) in Theorem \ref{linear pro}, 
is reduced to solving a linear programming problem. 
Therefore we may consider Theorem \ref{linear pro} as a linear programming bound in Euclidean space. 

Next we set $g_{ij}(r):=a_{ij}r^{2i+j}\;(a_{ij}\not=0)$ 
if $2i+j\le t$, 
and $g_{ij}(r):=0$ otherwise. 
Then we see by Lemma \ref{eujo} that the condition (\ref{to}) is equivalent to that 
$(X,w)$ is a Euclidean $t$-design. 
Therefore we see that a Euclidean $t$-design minimizes the potential energy 
for the functions 
$h(r,s,t)=g_0(r,s)+\sum_{i,j\ge 1} g_{ij}(r) g_{ij}(s)Q_j(t)$. 
The purpose of this paper is to introduce the following concept: 
\begin{defi}\label{seddef}
$(X,w)$ is called a \textit{strong Euclidean $t$-design} if 
the following condition holds:
\begin{equation}\label{sed}
\sum_{x,y \in X} w(x)w(y)(||x||||y||)^{i}Q_j(x'\cdot y')=0 \; ( 0\le\forall i\le t, \;1\le \forall j\le t).
\end{equation}
\end{defi}
strong Euclidean $t$-designs can be interpreted as the strongest designs among those minimizing the potential energy 
in the case when we take monomials as $g_{ij}$'s in Theorem \ref{linear pro}. 
For Euclidean designs, the Fisher type inequality is famous in 
algebraic combinatorics and numerical analysis. 
(cf. \cite{{Del},{M},{M2},{BBHS}}). 
Therefore it is natural to ask whether a Fisher type 
inequality holds for strong Euclidean designs, too. 
The following is the main result of this paper: 
\begin{theo}\label{fisher}
Let $(X,w)$ be a strong Euclidean $t$-design 
on $p$ concentric spheres. Assume that $p\ge e+1$. Then 
the following inequality holds: 
\[|X|\ge 
\begin{cases} 
\displaystyle (e+1)\left\{{d+e-1\choose e}+{d+e-2 \choose e-1}\right\} & \mbox{if $t=2e$},\\
\displaystyle 2(e+1){d+e-1 \choose e}& 
\begin{array}{l}
\mbox{if $(X,w)$ is antipodal and}\\
\mbox{$t=2e+1$,}\\
\end{array} 
\end{cases}\]
where we say $(X,w)$ is antipodal if $X=-X,\;w(x)=w(-x)$. 
\end{theo}
We give a proof of Theorem \ref{fisher} in Section \ref{fisher type}. 

Let $(X,w)$ be a strong Euclidean $t$-design on $p$ concentric spheres and assume $p\le t+1$. 
Then we see by Lemma \ref{spht} that each $X_i$ is a spherical $t$-design. 
On the other hand, it is easy to see that 
if $X_i$ is a spherical $t$-design for any $i$, 
then $X$ is a strong Euclidean $t$-design. 
Hence when $p\le t+1$, strong Euclidean $t$-designs 
are characterized by the property that each $X_i$ is a spherical $t$-design. 
So it is essential to consider the case when $p\geq t+2$. 
Moreover if a tight spherical $t$-design exists on $S^{d-1}$ 
(for the definition of tight spherical $t$-design, 
see \cite{DGS} or Theorem \ref{dgsfisher} in Section \ref{fisher type} of this paper), 
then putting a tight spherical $t$-design on each $e+1$ concentric spheres, 
we obtain an example attaining the lower bound in Theorem \ref{fisher}. 
\begin{defi}
A finite set in Euclidean space is called a {\it Euclidean $a$-code} 
if any distinct two points in the set are separated at least 
by $a$. 
%$X \subset \R^d$ is an {\it Euclidean $a$-code} if 
%$||x-y||^2\ge a^2$ $(x,y \in X,\ x\neq y)$ holds. 
\end{defi}
Finally in Section \ref{ECD},
as a corollary of Theorem \ref{linear pro}, 
we give a method to calculate an upper bound of the cardinality 
of a Euclidean $a$-code under the condition that the radii of concentric spheres on which the code lies are given. 
\section{Linear programming bound}\label{linearse}
First we give an equivalent condition for $(X,w)$ to be a Euclidean $t$-design: 
\begin{lem}[Neumaier-Seidel \cite{Neu}]\label{neujo}
Let $X$ be a finite set in $\R^d\backslash \{0\}$ 
and $w:X\rightarrow \R_{>0}$. 
Then the following are equivalent:
\begin{enumerate}
\item [{\rm (1)}]

$(X,w)$ is a Euclidean $t$-design. 

\item [{\rm (2)}] For any $1\le l \le t$,
$\phi\in \ha_l(\R^d)$ and $0\le j \le \left[\frac{t-l}{2}\right]$, 
\[
\sum_{x \in X} w(x)||x||^{2j}\phi(x)=0. 
\]
\end{enumerate}
\end{lem}
We define a non-degenerate inner product in the space $P_l(\R^d)$ 
as follows: 
for $f$, $g\in P_l(\R^d)$, 
\[\langle f,g \rangle :=\int_{S^{d-1}}f(x)g(x) d\sigma(x).\]
In a similar way to the above, define an inner product in the vector spaces $\ho_i(\R^d)$ 
and $\ha_i(\R^d)$. 
The following addition formula of the Gegenbauer polynomials 
is well-known: 
\begin{lem}[cf. \cite{DGS}]\label{kaho}
Let $\{\phi_{l,1},\ldots \phi_{l,h_l}\}$ be an 
orthonormal basis of $\ha_{l}(\R^d)$. 
Then for any $x$, $y\in S^{d-1}$, we have 
\[\sum_{i=1}^{h_l}\phi_{l,i}(x)\phi_{l,i}(y)=Q_l(x\cdot y).\]
\end{lem}
Also the following lemma is well-known: 
\begin{lem}[cf. \cite{S}]\label{seiteichi}
For any non-negative integer $i$ and any finite subset $X \subset S^{d-1}$, 
the matrix $\big(Q_i(x \cdot y)\big)_{x,y\in X}$ indexed by $X \times X$ is 
positive semi-definite. 
\end{lem}
We state an equivalent condition for $(X,w)$ to be a Euclidean $t$-design, in terms of the Gegenbauer polynomials. 
\begin{lem}\label{eujo}
$(X,w)$ is a Euclidean $t$-design if and only if 
for $1\le l \le t$ and $0\le j \le \left[\frac{t-l}{2}\right]$, 
the following equation holds: 
\[\sum_{x,y\in X} w(x)w(y)(||x||||y||)^{l+2j}Q_l\left(x'\cdot y'\right)=0.\]
\end{lem}
\begin{proof}
Let $\{\phi_{l,1},\ldots \phi_{l,h_l}\}$ be an 
orthonormal basis of $\ha_{l}(\R^d)$. Then 
\begin{eqnarray*}
\sum_{i=1}^{h_l}\left(\sum_{x\in X} w(x)||x||^{2j}\phi_{l,i}(x)\right)^2&=&\sum_{i=1}^{h_l}\sum_{x,y\in X} w(x)w(y)(||x||||y||)^{2j} \phi_{l,i}(x)\phi_{l,i}(y)\\
&=&\sum_{i=1}^{h_l}\sum_{x,y\in X} 
w(x)w(y)(||x||||y||)^{2j+l} \phi_{l,i}(x')\phi_{l,i}(y')\\
&=&\sum_{x,y\in X} w(x)w(y)(||x||||y||)^{2j+l}Q_l(x'\cdot y'). 
\end{eqnarray*}
The last equality follows from Lemma \ref{kaho}. 
By Lemma \ref{neujo}, 
$X$ is a Euclidean $t$-design if and only if 
$\sum_{i=1}^{h_l}\left(\sum_{x\in X} w(x)||x||^{2j}\phi_{l,i}(x)\right)^2=0$ 
holds for $1\le l \le t$ and $0\le j \le \left[\frac{t-l}{2}\right]$. 
Therefore the proof is completed. 
\end{proof}
\begin{proof}
[\it Proof of Theorem \ref{linear pro}]
By Definition \ref{pe}, 
\[PE_{f}(X)=\sum_{x,y \in X,\ x\neq y} w(x)w(y)f(||x||,||y||, x'\cdot y').\]
Since 
$f(r,s,t)\ge h(r,s,t)$ for any $r$, $s$ and $t$ and since 
$w(x)>0$ for any $x\in X$, we have 
\begin{eqnarray*}
PE_{f}(X) &\ge& \sum_{x \neq y \in X} w(x)w(y)h(||x||,||y||, x'\cdot y')\\
&=& \sum_{x,y \in X} w(x)w(y)h(||x||,||y||, x'\cdot y')-\sum_{x\in X}w(x)^2h(||x||,||x||,1)\\
&=& \sum_{x,y \in X} w(x)w(y)g_0(||x||,||y||)+\sum_{x,y \in X}
\sum_{i,j\ge 1}
w(x)w(y)g_{ij}(||x||)g_{ij}(||y||)Q_j(x' \cdot y')\\
&&-\sum_{x\in X}w(x)^2h(||x||,||x||,1).\\
\end{eqnarray*}
By Lemma \ref{seiteichi}, 
$\big(Q_i(x \cdot y)\big)$ is positive semi-definite. 
Hence for any $j\ge 1$, we have 
\[w(x)w(y)g_{ij}(||x||)g_{ij}(||y||)Q_j(x' \cdot y')\ge 0.\]
Namely 
\begin{eqnarray*}
PE_{f}(X)&\ge& \sum_{x,y \in X} w(x)w(y)g_0(||x||,||y||)-\sum_{x\in X}w(x)^2h(||x||,||x||,1).\\
\end{eqnarray*}
Therefore the inequality (\ref{linear}) holds. 
The condition to satisfy equality in (\ref{linear}) is clear. 
\end{proof}
Now by using Theorem \ref{linear pro}, it is easy to prove 
the generalized Sidelnikov inequality in Theorem \ref{gs}. 
First we quote the following well-known lemma: 
\begin{lem}[\cite{BB}, Lemma 3.4.3]\label{banho}
Let $t^l=\sum_{i=0}^l A_{l,i} Q_{l-i}(t)$ be 
the expansion in the Gegenbauer polynomials. 
Then 
\[A_{l,i}=
\begin{cases} 
\displaystyle \frac{l!!(d-2)!!}{i!!(d+2l-i-2)!!} & \mbox{if $i\equiv 0 \pmod{2}$},\\
\displaystyle 0 & \mbox{otherwise}.
\end{cases}\]
\end{lem}
\begin{proof}
[Proof of Theorem \ref{gs}]
In Theorem \ref{linear pro}, 
we set $f(r,s,t)=(rs)^{l}t^l$, 
\[
\left\{
\begin{array}{ll}
g_{0j}(r)=\sqrt{A_{l,l-j}}r^{l} &0\le j \le l\\
g_{ij}(r)=0 & i\not=0
\end{array}
\right.
\]
and $g_{0}(r,s)=g_{00}(r)g_{00}(s)$. Then by Lemma \ref{banho}, 
\begin{align*}
h(r,s,t)&= g_0(r,s)+\sum_{i,j\ge 1} g_{ij}(r)g_{ij}(s)Q_j(t)\\
&=\sum_{j=0}^l A_{l,l-j} (rs)^{l}Q_{j}(t)
=(rs)^{l}t^l=f(r,s,t).
\end{align*}
Therefore the conditions in Theorem \ref{linear pro} hold. Hence 
\begin{eqnarray*}
PE_{f}(X)&=&\sum_{x\not=y \in X} w(x)w(y)f(||x||,||y||,x'\cdot y')=\sum_{x\not=y \in X} w(x)w(y)(x\cdot y)^l\\
&\ge& \sum_{x,y \in X} w(x)w(y)g_0(||x||,||y||)-\sum_{x\in X}w(x)^2h(||x||,||x||,1)\\
&=& \sum_{x,y \in X} w(x)w(y)A_{l,l}(||x||||y||)^l-\sum_{x\in X} w(x)^2f(||x||,||x||,1)\\
&=& A_{l,l}\left(\sum_{x \in X} w(x)||x||^l\right)^2-\sum_{x \in X} w(x)^2 f(||x||,||x||,1).
\end{eqnarray*}
Since $A_{l,l}=A_l$, we have 
\[\sum_{x,y\in X} w(x)w(y)(x\cdot y)^l \ge A_l{\left(\sum_{x \in X}w(x)||x||^l\right)^2}_.\]
\end{proof}
\section{Fisher type inequality}\label{fisher type}
In this section, we give a Fisher type inequality for strong Euclidean designs. First 
we show the following lemma: 
\begin{lem}\label{spht}
Let $(X,w)$ be a strong Euclidean $t$-design on $p$ concentric spheres, 
and $w:X\rightarrow \R_{>0}$ be constant on each concentric sphere. 
Suppose $t+1\ge p$. Then each $X_i$ is 
a spherical $t$-design. 
\end{lem}
\begin{proof}
Let $(X,w)$ be a strong Euclidean $t$-design. Then 
by Definition \ref{seddef}, for $0\le i\le t$ and $1\le  j\le t$
\[\sum_{x,y \in X} w(x)w(y)(||x||||y||)^{i}Q_j(x'\cdot y')=0.\]
Let $\{\phi_{l,1},\ldots \phi_{l,h_l}\}$ be an 
orthonormal basis of $\ha_{l}(\R^d)$. Then for any $j$, 
\begin{eqnarray*}
\sum_{i=1}^{h_l}\left(\sum_{x\in X} w(x)||x||^j\phi_{l,i}(x)\right)^2&=&\sum_{i=1}^{h_l}\sum_{x,y\in X} w(x)w(y)(||x||||y||)^{j} \phi_{l,i}(x)\phi_{l,i}(y)\\
&=&\sum_{i=1}^{h_l}\sum_{x,y\in X} w(x)w(y)(||x||||y||)^{j+l} \phi_{l,i}(x')\phi_{l,i}(y')\\
&=&\sum_{x,y\in X} w(x)w(y)(||x||||y||)^{j+l}Q_l(x'\cdot y').
\end{eqnarray*}
The last equality follows from Lemma \ref{kaho}. 
Therefore $(X,w)$ is a strong Euclidean $t$-design 
if and only if the following equalities hold: 
for $1\le l \le t$, $\phi(x)\in \ha_{l}(\R^d)$ and 
$-l \le j \le t-l$, 
\begin{align}\label{sed1}
\sum_{x\in X}w(x)||x||^j\phi(x)&=\sum_{i=1}^p |r_i|^j \sum_{x\in X_i} w(x)\phi(x)=0. 
\end{align} 
Fix $\phi \in \ha_{l}(\R^d)$ and 
regard $\{\sum_{x\in X_i} w(x)\phi(x)\}_{i=1}^p$ as variables. Then 
the matrix coefficient of the linear system (\ref{sed1}) 
is 
\[\left( \begin{array}{cccc} r_1^{-l} & r_2^{-l} & \ldots & r_p^{-l}\\
r_1^{-l+1} & r_2^{-l+1} & \ldots & r_p^{-l+1}\\
\vdots & \vdots & \vdots & \vdots\\
r_1^{t-l} & r_2^{t-l} & \ldots & r_p^{t-l} \end{array} \right)_. \] 
When $t+1\ge p$, the rank of this matrix is $p$. 
Hence for any $1\le i \le p$ and $\phi \in \ha_{l}(\R^d)$, we have 
$$\sum_{x\in X_i} w(x)\phi(x)=0.$$ 
Now because of the fact that 
$w(x)$ is constant on each concentric sphere, 
each $X_i$ is a spherical $t$-design by Lemma \ref{sdec}. 
\end{proof}
The following theorem is the well-known Fisher type inequality for spherical designs: 
\begin{theo}[Delsarte-Goethals-Seidel \cite{DGS}]\label{dgsfisher}
Let $X \subset S^{d-1}$ be a spherical $t$-design. Then 
\begin{align}\label{ineq:Fisher}
|X|\ge 
\begin{cases} 
\displaystyle {d+e-1\choose e}+{d+e-2 \choose e-1} & \mbox{if $t=2e$},\\
\displaystyle 2{d+e-1 \choose e}& \mbox{if $t=2e+1$}.
\end{cases}
\end{align}
\end{theo}
A spherical $t$-design $X$ is {\it tight} if 
equality holds in (\ref{ineq:Fisher}). 
%Therefore, 
By Lemma \ref{spht} and Theorem \ref{dgsfisher}, 
we obtain the following corollary: 
\begin{cor}
Let $(X,w)$ be a strong Euclidean $t$-design on $p$ concentric spheres. 
Suppose that $w:X\rightarrow \R_{>0}$ be constant on each concentric sphere and that $t+1\ge p$. Then 
\[|X|\ge 
\begin{cases} 
\displaystyle p\left\{{d+e-1\choose e}+{d+e-2 \choose e-1}\right\} & \mbox{if $t=2e$},\\
\displaystyle 2p{d+e-1 \choose e}& \mbox{if $t=2e+1$}.
\end{cases}\]
\end{cor}

In the sequel, suppose that $p$ is sufficiently large comparing to $t$. 
Our proof below follows Delsarte-Seidel {\cite{Del}}. 
For a subspace $P$ of $\Pol(\R^d)$, put $||x||^{j}P:=\{||x||^jf(x) \mid f\in P\}$. 
We set \[\Pol'(\R^d):=\Pol(\R^d)+||x||\Pol(\R^d).\] 
We remark that the sum of the right hand side is a direct sum. 
It is because, if there exist nonzeros $f$,  $g\in \Pol(\R^d)$ 
such that $f+||x||g=0$, then 
$f^2=||x||^2g^2$. Because $||x||^2=x_1^2+\cdots+x_d^2$ is irreducible in 
$\Pol(\R^d)$, we have a contradiction since 
the parities of $||x||^2=x_1^2+\cdots+x_d^2$ 
in the left and right hand side are different. 
Set 
\begin{align*}
\Pol'_j(\R_d)&:=\Pol_j(\R_d)+||x||\Pol_{j-1}(\R_d),\\
\ho_j'(\R^d)&:=\ho_j(\R^d)+||x||\ho_{j-1}(\R^d).
\end{align*}
Then 
\[\Pol'_j(\R^d)=\bigoplus_{i=0}^j \ho'_j(\R^d).\]
Generally for $T \subset \R^d$, we denote by $\ho_l(T)$ (resp.\ $\ha_l(T)$) 
the vector space of elements of $\ho_l(\R^d)$ (resp.\ $\ha_l(\R^d)$) 
which are restricted on $T$. 
For example we write $\Pol(T)=\{f|_{T} \mid f\in \Pol(\R^d)\}$, 
where $f|_T$ denotes a restricted function on $T$ for $f$. 
\begin{lem}\label{j}
\[\ho'_j(RS)=\ho_j(RS)\oplus (||x||\ho_{j-1})(RS).\]
\end{lem}
\begin{proof}
Take any $f\in \ho_j(RS)$ and $g\in \ho_{j-1}(RS)$ such that 
$f=||x||g$, 
then we have $f^2(x)=||x||^2g^2(x)$ as polynomials. 
Since $||x||^2=x_1^2+x_2^2+\cdots+x_d^2$ is an irreducible element of 
the polynomial ring, 
checking the parities of $||x||^2=x_1^2+x_2^2+\cdots+x_d^2$ in 
the left and right hand side, we have $f=g=0$. 
Therefore the sum of the right hand side is a direct sum. 
\end{proof}
\begin{lem}
Suppose that $RS$ consists of $p$ concentric spheres. Then we have 
\begin{equation}\label{bun}
\ho'_j(RS)\subset \sum_{i=1}^p \ho'_{j+i}(RS).
\end{equation}
\end{lem}
\begin{proof}
For $f\in \ho_j'(\R^d)$, we have the following identity on $RS$: for 
$y\in RS$ 
\begin{align}\label{eqn:exp}
f(y)\prod_{r\in RX}(r-||y||)=0.
\end{align}
Expanding (\ref{eqn:exp}), we see that 
$f(y)$ is written as a linear combination 
with respect to $||y||^i f(y)\;(i=1,2,\ldots,p)$, where $||y||^if(y) \in \ho_{j+i}'(RS)$. 
\end{proof}
\begin{lem}\label{chokuwa}
Suppose that $RS$ consists of $p$ concentric spheres. Then we have 
\[\Pol'_j(RS)=\bigoplus_{i=0}^{p-1}\ho'_{j-i}(RS)\]
\end{lem}
\begin{proof}
By (\ref{bun}), we have 
\begin{equation}\label{wa}
\Pol'_j(RS)=\sum_{i=0}^j\ho'_{i}(RS)=\sum_{i=0}^{p-1} \ho'_{j-i}(RS).
\end{equation}
Therefore it is enough to show that 
the sum of the right hand side is a direct sum. 
First we show that for the restriction homomorphism 
$\phi:\Pol'_j(\R^d)\rightarrow \Pol'_j(RS)$, 
\begin{equation}\label{ker}
\Ker\;\phi=\Pol'_{j-p}(\R^d)\prod_{r\in RX}(r-||x||). 
\end{equation}
Clearly we have 
$$\Ker\;\phi \supset \Pol'_{j-p}(\R^d)\prod_{r\in RX}(r-||x||).$$ 
Conversely, 
take $f+||x||g \in \Ker\;\phi,\;(f\in \Pol_j(\R^d),\;g\in \Pol_{j-1}(\R^d))$. 
For $r_1 \in RX$, $$f(x)+||x||g(x)=f(x)+r_1g(x)-(r_1-||x||)g(x).$$ 
Hence $f(x)+r_1g(x)$ is zero on $r_1S^{d-1}$. 
By Hilbert's Nullstellensatz, there exists some $h(x)\in \Pol_{j-2}(\R^d)$ 
such that $f(x)+r_1g(x)=(r_1^2-||x||^2)h(x)$. 
Therefore we have 
$$f(x)+||x||g(x)=(r_1-||x||)\left\{(r_1+||x||)h(x)-g(x)\right\}.$$ 
Similarly, replacing 
$f(x)+||x||g(x)$ and $r_1$ by $(r_1+||x||)h(x)-g(x)$ and $r_2$, respectively,
we see that there exists $q(x)\in \Pol'_{j-2}(\R^d)$ such that 
$$f(x)+||x||g(x)=(r_1-||x||)(r_2-||x||)q(x).$$ 
Recursively we see that there exists $r(x)\in \Pol'_{j-p}(\R^d)$ such that 
$$f(x)+||x||g(x)=\prod_{r\in RX}(r-||x||)r(x).$$ 
Therefore we have $\Ker\phi \subset \Pol'_{j-p}(\R^d)\prod_{r\in RX}(r-||x||)$. 

By (\ref{ker}), 
\[\dim \Pol'_j(RS)=\dim \Pol'_j(\R^d)-\dim \Pol'_{j-p}(\R^d).\]
Using $\ho'_i(\R^d)\simeq \ho'_i(RS)$, we have 
\[\dim \Pol'_j(RS)=\sum_{i=0}^{p-1} \dim \ho'_{j-i}(\R^d)=\sum_{i=0}^{p-1} \dim \ho'_{j-i}(RS).\]
This implies that the sum of the right hand side in (\ref{wa}) is 
a direct sum. 
\end{proof}
\begin{proof}
[Proof of Theorem \ref{fisher}]
The following decomposition is well-known (cf. \cite{Neu}): 
\[\ho_i(\R^d)=\bigoplus_{j=0}^{[\frac{i}{2}]}||x||^{2j}\ha_{i-2j}(\R^d).\]

We set 
\[P_t:=\bigoplus_{i=0}^t \bigoplus_{j=-i}^{t-i} ||x||^j \ha_i(\R^d)=\sum_{i=0}^t\sum_{j=-i}^{t-i}||x||^j\ho_{i}(\R^d).\] 
Then for any $f(x)= \sum_{i=0}^t\sum_{j=-i}^{t-i}||x||^jf_i(x)\in P_t$ with $f_i \in \ha_i(\R^d)$, 
\begin{eqnarray*}
&&\sum_{k=1}^p \frac{w(X_k)}{|S_k|}\int_{S_k} f(x)d\sigma_k(x)=\sum_{k=1}^p\sum_{j=-i}^{t-i}r_k^j \frac{w(X_k)}{|S_k|}\sum_{i=0}^t\int_{S_k} f_i(x)d\sigma_k(x)
\\
&=& \sum_{k=1}^p\sum_{j=-i}^{t-i}r_k^j \frac{w(X_k)}{|S_k|}\int_{S_k} f_0(x)d\sigma_k(x)
=\sum_{k=1}^p\sum_{j=-i}^{t-i}r_k^j w(X_k)f_0(x)\\
&=&\sum_{k=1}^p\sum_{j=-i}^{t-i}\sum_{x\in X_k}w(x)||x||^jf_0(x)
= \sum_{j=-i}^{t-i}\sum_{x \in X} w(x)||x||^jf_0(x). 
\end{eqnarray*} 
Let $(X,w)$ be a strong Euclidean $t$-design. Then 
by the equivalent condition (\ref{sed1}) for $(X,w)$ to be a strong Euclidean $t$-design, 
we have 
\[\sum_{j=-i}^{t-i}\sum_{x \in X} w(x)||x||^jf_0(x)=\sum_{i=0}^t \sum_{j=-i}^{t-i} \sum_{x\in X} w(x)||x||^jf_i(x)= \sum_{x\in X} w(x)f(x).\]
Therefore for any $f\in P_t$, 
\begin{equation}\label{nai}
\sum_{k=1}^p \frac{w(X_k)}{|S_k|}\int_{S_k} f(x)d\sigma_k(x)=\sum_{x\in X} w(x)f(x). 
\end{equation}
Suppose $t=2e$. Then we have 
$P_t=P_eP_e=\langle f\cdot g\mid f,g \in P_e\rangle$, 
where $\langle f\cdot g\mid f,g \in P_e\rangle$ is 
the vector space expanded by $f\cdot g\; (f,g \in P_e)$. 
 
We define the non-degenerate inner products $[\cdot,\cdot]$ and $\langle \cdot,\cdot\rangle_{RS}$ 
on $P_e(X)$ and $P_e(RS)$, respectively, as follows: for $f$, $g \in P_e$,
\begin{align}
[f,g]&:=\sum_{x\in X} w(x) f(x)g(x), \label{naiseki1} \\
\langle f,g \rangle_{RS}&:=\sum_{i=1}^p \frac{w(X_i)}{|S_i|}\int_{S_i} f(x)g(x)d\sigma_i(x). \label{naiseki2}
\end{align}
Then (\ref{nai}) is equivalent to that, for any $f$, $g \in P_e$,
\[[f,g]=\langle f,g \rangle_{RS}.\]
This implies that the restriction map 
$\rho_e:P_e(RS) \rightarrow P_e(X)$ is an injective homomorphism. 
Hence $|X|$ is bounded below by $\dim P_e(RS)$. 

Set $T_e:=(||x||^eP_e)(RS)=\sum_{i=0}^t\sum_{j=-i}^{t-i}\left(||x||^{e+j}\ho_{i}\right)(RS)$. By the fact that $\dim P_e(RS)=\dim T_e$, 
it is enough to calculate $\dim T_e$. 
Generally we have the following: 
\begin{equation}\label{2age}
(||x||^2\ho_{i-2})(RS) \subset \ho_i(RS).
\end{equation}
Therefore 
\[T_e=\sum_{i=1}^{e+1} (||x||^i\ho_{e-1})(RS)+\sum_{i=0}^{e}(||x||^i\ho_{e})(RS).\]
By the assumption $p \ge e+1$, 
this sum is a direct sum by Lemma \ref{j} and \ref{chokuwa}. So we have 
\[\dim T_e=(e+1)\left\{{d+e-1\choose e}+{d+e-2\choose e-1}\right\}_.\]
Next we suppose that $(X,w)$ is antipodal and that $t=2e+1$. 
Then set 
\[P_{2e}'(\R^d)=\bigoplus_{i=0}^e\bigoplus_{j=-2i}^{2e+1-2i}||x||^j \ha_{2i}(\R^d)=\sum_{i=0}^e\sum_{j=-2i}^{2e+1-2i}||x||^j \ho_{2i}(\R^d).\] 
We assume that $X$ is a disjoint union of $Y$ and $-Y$. 
Then in a similar way to the case when $t=2e$, we see that 
$(X,w)$ is a strong Euclidean $(2e+1)$-design if and only if 
for any $f\in P'_{2e}$ the following holds: 
\begin{equation}\label{nai2}
\sum_{i=1}^p \frac{w(X_i)}{|S_i|}\int_{S_i} f(x)d\sigma_i(x)=\sum_{y\in Y} w(y) f(y). 
\end{equation}
Set 
\[P_{e}''(\R^d)=\sum_{k=0}^{[\frac{e}{2}]}\sum_{j=-e+2k}^{2k}||x||^j \ho_{e-2k}(\R^d).\]
Then we have $P_{2e}'(\R^d) \supset P_e''(\R^d)\cdot P''_e(\R^d)$. 
Therefore when we also define the non-degenerate inner products in the space 
$P_e''(Y)$ and $P_e''(RS)$ in the same way as (\ref{naiseki1}) and (\ref{naiseki2}),  
we see by (\ref{nai2}) that, for any $f$, $g\in P''_e$ 
\[[f,g]=\langle f,g \rangle_{RS}.\]
Therefore we see that the restriction map $\rho'_e:P''_e(RS) \rightarrow P''_e(Y)$ 
is injective. 
Hence, $|Y|$ is bounded below by $\dim P''_e(RS)$. 
In particular, $|X|=2|Y|\ge 2\dim P''_e(RS)$. 
Set $T'_e:=(||x||^eP''_e)(RS)$. 
Then from the fact that $\dim P''_e(RS)=\dim T'_e$ and from (\ref{2age}), we have 
\begin{align}\label{eqn:1.3_anti}
T_e'=\ho_e(RS)+||x||\ho_e(RS)+\cdots+||x||^e\ho_e(RS).
\end{align}
When $p\ge e+1$, by Lemma \ref{chokuwa}, the sum of the 
right hand side in (\ref{eqn:1.3_anti}) 
is a direct sum. 
Therefore $\dim T_e'=(e+1){d+e-1 \choose e}$. 
\end{proof}
\section{Bounds for Euclidean $a$-codes and Euclidean designs}\label{ECD}
In this section, we give a method 
to obtain a bound of the cardinality of Euclidean $a$-codes and 
Euclidean designs. 
\begin{theo}\label{upper}
Let $X \subset \R^d\;(|X|<\infty)$ be a Euclidean $a$-code, 
$g_0:\R_{>0}^2\rightarrow \R$ , $g_{ij}:\R_{>0}\rightarrow \R \;(i,j \ge 1)$ and $h(r,s,t)=g_0(r,s)+\sum_{i,j\ge 1}g_{ij}(r)g_{ij}(s)Q_j(t)$. 
Assume that 
\[h(r,s,t)\le 0, \;(r^2+s^2-2rst\ge a^2,\; r,s\in RX).\]
Then we have the following inequality: 
\begin{equation}\label{a}
\sum_{x,y \in X} g_{0}(||x||,||y||)\le \sum_{x\in X} h(||x||,||x||,1).
\end{equation}
\end{theo}
\begin{proof}
Set 
\[f(r,s,t)=
\begin{cases} 
+\infty & \mbox{if $r^2+s^2-2rst< a^2$},\\
0 & \mbox{otherwise}.
\end{cases}\]
Since  $h(r,s,t)\le 0, \;(r^2+s^2-2rst\ge a^2,\; r,s\in RX)$ by the assumption,
we have 
$f(r,s,t)\ge h(r,s,t),\;(r,s \in RX, t \in [-1,1])$. 
Therefore, setting $w\equiv 1$, we have the following inequality by Theorem \ref{linear pro}, 
\begin{equation}\label{katei}
PE_f(X,w)\ge \sum_{x,y \in X} g_{0}(||x||,||y||) - \sum_{x\in X} h(||x||,||x||,1).
\end{equation}
If $\sum_{x,y \in X} g_{0}(||x||,||y||) > \sum_{x\in X} h(||x||,||x||,1)$, 
then the right hand side of (\ref{katei}) is positive, and so $PE_f(X,w)= +\infty$. 
Hence, there exist $x \neq y \in X$ such that $||x-y||^2\le a^2$. 
This contradicts to that $X$ is a Euclidean $a$-code. 
\end{proof}
\begin{cor}\label{sei}
Let $X \subset \R^d$ be a Euclidean $a$-code, 
$RX=\{r_1,\ldots,r_p\}$ and $X_i=X\cap r_i S^{d-1}$. 
Set $f_{ij}(r,s,t)=(rs)^iQ_j(t),\;(i, j\ge 0)$, and $h(r,s,t)=\sum_{i,j}a_{ij}f_{i,j}(r,s,t)$, $a_{ij}\ge 0 \;(\forall j \ge 1)$. 
Assume that $h(r,s,t)\le 0, \;(r^2+s^2-2rst\ge a^2,\; r,s\in RX)$. 
Then the following inequality holds: 
\begin{equation}\label{a1}
\sum_{i}a_{i0}\left(\sum_{k=1}^p r_k^i|X_k|\right)^2\le \sum_{k=1}^p h(r_k,r_k,1)|X_k|.
\end{equation}
\end{cor}
\begin{proof} 
Set $g_{ij}(r)=\sqrt{a_{ij}}r^i\;(j\ge 1)$ and $w \equiv 1$, and  use Theorem \ref{upper}. 
\end{proof}
In (\ref{a1}), the left hand side is of degree two 
and the right hand side is of degree one with respect to $|X_i|$. 
So (\ref{a1}) gives an upper bound of $|X_i|$ if we can find a good function. 
\begin{remark}
Let $X \subset \R^d$ be a Euclidean $a$-code and 
$X':=\{x'\mid x \in X\}$. 
If $|X'|<|X|$, namely if there exist $x\neq y \in X$ such that 
$x'=y'$, then $X_i$ and $X_j$ are separated at least by $a$ where $x\in X_i$ and $y \in X_j$. 
Hence, the condition of a Euclidean $a$-code does not give any restriction between $X_i$  and $X_j$.
Therefore, it is enough to consider the case $|X'|=|X|$. 
For $x,y \in X, (||x||=r, ||y||=s)$, we have $||x-y||^2=r^2+s^2-2rst \ge a^2$. 
For $r,s \in RX$, set
\[z_{rs}:=\frac{r^2+s^2-a^2}{2rs},\;z:=\max\{z_{rs}\mid r,s\in RX \}.\]
Then since $X$ is a Euclidean $a$-code, we have $x'\cdot y'\le z\ (\forall x,y \in X)$. 
So $X'$ is a spherical $z$-code. 
The linear programming bound of the usual Delsarte method for spherical $z$-codes is a method as 
 seeking a polynomial $f(t)=\sum_i a_i Q_i(t)$ such that $a_i\ge 0\ (\forall i \ge 1)$, 
$a_0>0$ and $f(t)\le 0 \;(\forall t \in [-1,z])$ (refer to \cite{Del1}). 
In Corollary \ref{sei}, if $a_{ij}=0\ (\forall i\ge 1)$, 
then (\ref{a1}) is the same to the linear programming bound of the Delsarte method. 
Since functions $f(t)$ or bounds appearing in the Delsarte bounds are particular cases in ones of Corollary \ref{sei}. So 
there is a possibility to improve the bound obtained by using the Delsarte method directly to Euclidean $a$-codes as above. 
\end{remark}
Now we give a linear programming bound on the 
cardinality of a Euclidean $t$-design. 
Set 
$A(X_{(r)},X_{(s)}):=\{x'\cdot y'\mid x\in X_{(r)}, y\in X_{(s)}, x'\not=y'\}$.
\begin{theo}\label{dlpb}
Let $(X,w)$ be a Euclidean $t$-design. 
Suppose that $w(x)$ is constant on each concentric sphere and 
denote by $w(||x||):=w(x)$. Set 
$$I:=\{(i,j)\in \Z^2_{\ge 0} \mid 0\le i \le t,\;
\mbox{there exists $k \in \Z_{\ge 0}$ such that $i=2k+j$ or $j=0$} \}.$$ 
Assume that for any 
$(i,j) \not\in I$, $a_{ij}\le 0$ and 
$f(r,s,t)=\sum_{i,j\ge 0}a_{ij} f_{ij}(r,s,t)$. 
Moreover, assume that 
$f(r,s,t) \ge 0 \;(\forall r,s \in RX,\;\forall t \in [-1,1])$. 
Then, we have the following inequality: 
\begin{align}\label{ineq:4.2}
\sum_{i\ge 0}a_{i0}\left(\sum_{r\in RX}w(r)r^i|X_{(r)}|\right)^2\ge \sum_{r,s\in R(X)} w(r)w(s)f(r,s,1)d_{r,s,1},
\end{align}
where 
\[d_{r,s,t}:=\sharp\{(x,y)\in X^2\mid x\in X_{(r)}, y\in X_{(s)}, x'\cdot y'=t\}.\]
\end{theo}
\begin{proof}
We estimate the following value: 
\begin{equation}\label{hyo}
\sum_{i,j} a_{ij} \sum_{x,y \in X} w(x)w(y)f_{ij}(||x||,||y||,x' \cdot y'). 
\end{equation}
Since $(X,w)$ is a Euclidean $t$-design, 
\begin{eqnarray*}
(\ref{hyo}) \;&=&\sum_{i\ge0}a_{i0}\sum_{x,y\in X} w(x)w(y)(||x||||y||)^i\\
&+&\sum_{(i,j)\not\in I} a_{ij}\sum_{x,y\in X}w(x)w(y)(||x||||y||)^iQ_j(x'\cdot y'). 
\end{eqnarray*}
By the assumption, 
$a_{ij}\le 0\;(\forall (i,j)\not\in I)$ and by Lemma \ref{seiteichi}, 
$$\sum_{x,y\in X}w(x)w(y)(||x||||y||)^iQ_j(x'\cdot y')\ge 0.$$ 
Hence, 
\[(\ref{hyo})\le \sum_{i\ge0}a_{i0}\sum_{x,y\in X} w(x)w(y)(||x||||y||)^i= \sum_{i\ge 0}a_{i0}{\left(\sum_{r\in RX}w(r)r^i|X_{(r)}|\right)^2}_.\]
On the other hand, 
by $f(r,s,t) \ge 0 \;(\forall r,s \in RX,\;\forall t \in [-1,1])$, we have 
\begin{eqnarray*}
(\ref{hyo})&=&\sum_{x,y \in X}w(x)w(y) \sum_{i,j} a_{ij} f_{ij}(||x||,||y||,x' \cdot y')\\
&=& \sum_{x,y \in X} w(x)w(y)f(||x||,||y||,x'\cdot y')\\
&\ge& \sum_{r,s\in R(X)} w(r)w(s)f(r,s,1)d_{r,s,1}
\end{eqnarray*}
\end{proof}
In (\ref{ineq:4.2}), the left hand side is of degree two 
and the right hand side is of degree one with respect to $|X_i|$. 
Therefore, (\ref{ineq:4.2}) give a good lower bound of $|X_i|$ 
for a suitable $f$ if we could find a good function $f$.

\end{document}